\documentclass[10pt,a4paper]{article}
\usepackage[latin1]{inputenc}
\usepackage{anysize}
\usepackage{graphicx}
\marginsize{2.5cm}{2.5cm}{2cm}{1cm}
\usepackage{amssymb,amsmath,amscd,amsfonts,amsthm,mathrsfs}
\newcommand{\cc}[1]{\mathcal{#1}}

\numberwithin{equation}{section}
\newtheorem{theorem}{Theorem}%{\bf}{\it }
\newtheorem{lemma}[theorem]{Lemma}%{\bf}{\it }
\newtheorem{conjecture}[theorem]{Conjecture}%{\bf}{\it }
\usepackage{url}
\begin{document}
\title{Decomposing the cube into paths}
\author{Joshua Erde \thanks{DPMMS, University of Cambridge}}
\maketitle
\begin{abstract}
We consider the question of when the $n$-dimensional hypercube can be decomposed into paths of length $k$. Mollard and Ramras \cite{MR2013} noted that for odd $n$ it is necessary that $k$ divides $n2^{n-1}$ and that $k\leq n$. Later, Anick and Ramras \cite{AR2013} showed that these two conditions are also sufficient for odd $n \leq 2^{32}$ and conjectured that this was true for all odd $n$. In this note we prove the conjecture.
\end{abstract}
\section{Introduction}
The \emph{$n$-dimensional hypercube} $\cc{Q}_n$ is a graph with vertex set $V = \{0,1\}^n$ and edge set $E=\{(x,y)\,:\, ||x-y||_1 = 1\}$. Problems of decomposing the hypercube into edge disjoint copies of smaller graphs have been considered by several authors, such as decomposing $\cc{Q}_n$ into trees \cite{F1990} \cite{WW2012}, into Hamiltonian cycles and matchings \cite{ABS1990} or into stars, $K_{1,r}$ for $r<n$ \cite{BEEH2001}. Mollard and Ramras \cite{MR2013}, motivated by applications in parallel processing (see \cite{L1992}), considered the problem of decomposing the hypercube into paths. A path of length $k$ is a sequence of distinct vertices $x_1,x_2, \ldots, x_{k+1}$ such that for all $1 \leq i \leq k$ $(x_i,x_{i+1}) \in E(\cc{Q}_n)$. Mollar and Ramras \cite{MR2013} noted that if $n$ is odd, and we wish to decompose $\cc{Q}_n$ into paths of length $k$, there are two simple necessary conditions that $k$ must satisfy. Firstly, since $|E(\cc{Q}_n)|=n2^{n-1}$ we must have that $k$ divides $n2^{n-1}$, which we write as $k\,|\,n2^{n-1}$. Secondly, since $\cc{Q}_n$ is $n$-regular, and $n$ is odd, each vertex must be the endpoint of at least one of the paths, and so we must have at least $2^{n-1}$ paths (since each path has $2$ endpoints). Therefore we must also have that $k \leq n$. Anick and Ramras \cite{AR2013} conjectuted:
\begin{conjecture}[\cite{AR2013}]
Let $n$ be odd and $k$ such that $k \,|\, n$ and $k \leq n$. Then $\cc{Q}_n$ can be decomposed into paths of length $k$.
\end{conjecture}

They showed that the conjecture holds for $n < 2^{32}$. The result of this note is to show that the conjecture holds for all $n$.

\begin{theorem}\label{t:main}
Let $n$ be odd and $k$ such that $k\,|\,n2^{n-1}$ and $k \leq n$. Then $\cc{Q}_n$ can be decomposed into paths of length $k$.
\end{theorem}

In the next section we provide a proof of Theorem \ref{t:main} and in the final section we briefly discuss what can be said about decomposing $\cc{Q}_n$ into paths of length $k$ for even $n$.

\section{Proof of Theorem \ref{t:main}}
A \emph{walk of length $k$} is a sequence of vertices $x_1,x_2, \ldots , x_{k+1}$, not necessarily distinct, such that for all $1 \leq i \leq k$ $(x_i,x_{i+1}) \in E(\cc{Q}_n)$. We will often define walks and paths by describing their edge sets. We denote by \emph{even vertices} the set of vertices $(q_1,q_2, \ldots , q_n) \in \cc{Q}_n$ such that $|\{i \,:\,q_i=1\}|$ is even, and similarly \emph{odd vertices}. It is apparent that $\cc{Q}_n$ is a bipartite graph, with the classes being the even and the odd vertices. Two vertices $x,y \in \cc{Q}_n$ are \emph{antipodal} if $||x-y||_1 = n$, and we call a path of length $n$ between two antipodal points an \emph{antipodal path}.\\

\begin{lemma}\label{l:antipodal}
For any $n$ $\cc{Q}_n$ can be decomposed into antipodal paths of length $n$.
\end{lemma}
\proof
Given a vertex $(q_1,q_2, \ldots , q_n) \in \cc{Q}_n$ there is a natural antipodal path to consider, that is 

\begin{align*}
\{\big(&(q_1,q_2, \ldots , q_n),(q_1 + 1 ,q_2, \ldots , q_n) \big), \big((q_1 + 1 ,q_2, \ldots , q_n),(q_1 + 1,q_2 + 1, q_3, \ldots , q_n)\big) , \ldots, \\
\big(&(q_1 +1 ,q_2 + 1, \ldots , q_{n-1}+1 , q_n),(q_1 +1 ,q_2 + 1, \ldots,q_{n-1}+1 , q_n + 1)  \big)\}
\end{align*}

where addition is taken modulo $2$. If we only take the paths beginning at even vertices then we cover each edge exactly once. Indeed if an edge $\left((p_1,p_2, \ldots , p_i, \ldots , p_n),(p_1,p_2, \ldots, p_i +1 , \ldots , p_n)\right)$ is in two of these paths, then we must have that $(p_1,p_2, \ldots , p_i, \ldots , p_n)$ is the $i$th vertex in one path and $(p_1,p_2, \ldots, p_i +1 , \ldots , p_n)$ is the $i$th vertex in the other. However this would imply that the number of $1$s in each vector has the same parity, a contradiction. Hence each edge is covered at most once and since there are $2^{n-1}$ even vertices, and each path has length $n$, we have covered $n2^{n-1} = |E(\cc{Q}_n)|$ edges. \qed \\
\ \\
Since we can decompose $\cc{Q}_n$ into paths of length $n$ it is also clear that we can decompose $\cc{Q}_n$ into paths of length $t$ for all $t\,|\,n$ by subdividing these antipodal paths in the natural way. In fact this simple observation achieves more if we consider the structure these paths induce on $\cc{Q}_n$. \\

\begin{lemma}\label{l:decomp}
For any $n$ let $t$ be such that $t$ is odd and $t\,|\,n$. If $\cc{Q}_{\frac{n}{t}}$ can be decomposed into paths of length $s$ then $\cc{Q}_n$ can be decomposed into paths of length $ts$.
\end{lemma}
\proof
Let us consider the antipodal paths on $\cc{Q}_n$ constructed in Lemma \ref{l:antipodal}. Suppose we split each of the paths into $\frac{n}{t}$ paths of length $t$. We define a graph $G$ on $\{0,1\}^n$ by joining two vertices if there is a path between them, that is, if one of the paths of length $t$ starts at one of the vertices and ends at the other. We claim that $G$ is just a disjoint union of copies of $\cc{Q}_{\frac{n}{t}}$. Indeed given a point $(q_1,q_2, \ldots , q_n)$ it is adjacent to the points
\begin{align*}
&(q_1+1,q_2+1, \ldots , q_t+1, q_{t+1}, \ldots, q_n), \\
&(q_1,q_2, \ldots , q_t, q_{t+1} +1, \ldots q_{2t}+1, q_{2t+1}, \ldots, q_n), \\
&\ldots \\
&(q_1,q_2, \ldots ,q_{n-t}, q_{n-t+1}+1, \ldots, q_n+1).
\end{align*}

So if we divide $\{0,1\}^n$ into equivalence classes under the relation $(q_1, q_2, \ldots , q_n) \sim (p_1, p_2, \ldots , p_n)$ if $(q_1 - p_1, q_2 - p_2 , \ldots , q_n - p_n) \in \{(0,0, \ldots, 0) , (1,1, \ldots , 1) \} ^{\frac{n}{t}}$ (where $(0,0, \ldots, 0)$ and $(1,1, \ldots , 1)$ are of length $t$), we see that $G$ restricted to each equivalence class is isomorphic to $\cc{Q}_{\frac{n}{t}}$, and each edge in $G$ is inside one equivalence class. \\
\ \\
We use the decomposition of $\cc{Q}_{\frac{n}{t}}$ into paths of length $s$ to decompose $G$ into paths of length $s$, and see that, when considered in $\cc{Q}_n$, a path of length $s$ in $G$ is a walk of length $ts$. More precisely if we have a path $\{(x_1,x_2) , (x_2,x_3) , \ldots, (x_{s-1}, x_s) \}$ in $G$ we know that each edge $(x_i,x_j)$ corresponds to some path of length $t$ in $\cc{Q}_n$, 

$$P_t^{i,j} = \{(x_i, x_{\{i,j\}_2}) , (x_{\{i,j\}_2}, x_{\{i,j\}_3}), \ldots ,(x_{\{i,j\}_{t-1}}, x_j) \}.$$ 

So we have that 

\begin{align*}
W = \{ &(x_1,x_{\{1,2\}_2}), \ldots , (x_{\{1,2\}_{t-1}}, x_2), (x_2,x_{\{2,3\}_2}), \ldots , \\
		&(x_{\{2,3\}_{t-1}}, x_3), (x_3,x_{\{3,4\}_2}) ,\ldots , (x_{\{s-1,s\}_{t-1}}, x_s) \}
\end{align*}			

is a walk in $\cc{Q}_n$. It remains to check that there are no repeated vertices in $W$.\\
\ \\
Since the decomposition of $\cc{Q}_{\frac{n}{t}}$ was into paths we know that $x_1,x_2, \ldots x_t$ are distinct and also we know the form that each $P_t^{i,j}$ takes. Given an interior point to a path, say $x_{\{i,j\}_{l}}$, we know that it agrees with $x_i$ and $x_j$ except in some subset of a block of $t$ consecutive co-ordinates (specifically differing in the first $l$, or the last $t-l$ of one of those). Hence given $x_{\{i,j\}_{l}}$ and the equivalence class of vertices we know $x_i$ and $x_j$, and so $x_{\{i,j\}_{l}}$ is interior to only one $P_t^{i,j}$. Since the points in each $P_t^{i,j}$ are distinct, and the interior points are not in the same equivalence class as the endpoints, we have that $W$ has no repeated vertices, and so is a path.\qed\\
\ \\

It follows by Lemma \ref{l:decomp} that we only need to consider the case of decomposing $\cc{Q}_n$ into paths of length $2^r$ for $2^r<n$. To prove this case we will need the following lemma.

\begin{lemma}\label{l:split}
If $\cc{Q}_i$ and $\cc{Q}_j$ can be decomposed into paths of length $k$ then so can $\cc{Q}_{i+j}$.
\end{lemma}
\proof
For each vertex $x \in \cc{Q}_i$ the subgraph of $\cc{Q}_{i+j}$ on the set of vertices $(q_1,q_2, \ldots, q_{i+j})$ such that $(q_{1},q_{2}, \ldots, q_{i}) = x$ is isomorphic to $\cc{Q}_j$, and so we can decompose each of these, disjoint, subgraphs, by using the decomposition of $\cc{Q}_j$. Similarly for each vertex $y \in \cc{Q}_j$ the subgraph of $\cc{Q}_n$ on the set of vertices $(q_1,q_2, \ldots, q_{i+j})$ such that $(q_{i+1},q_{i+2}, \ldots, q_{i+j}) = y$ is isomorphic to $\cc{Q}_i$ and so we can decompose these subgraphs by using the decomposition of $\cc{Q}_i$. Note that each edge is in exactly one of these subgraphs, since any edge is between two vertices which differ in exactly one co-ordinate, which is either in the first $i$, or the last $j$. \qed \\
\ \\
We will also need the following folklore result, for a proof see e.g. \cite{ABS1990}.
\begin{lemma}\label{l:Ham}
Let $n$ be even. Then $\cc{Q}_n$ can be decomposed into edge disjoint Hamiltonian cycles.
\end{lemma}

\begin{figure}[!ht]
\centering
\includegraphics[scale=0.2]{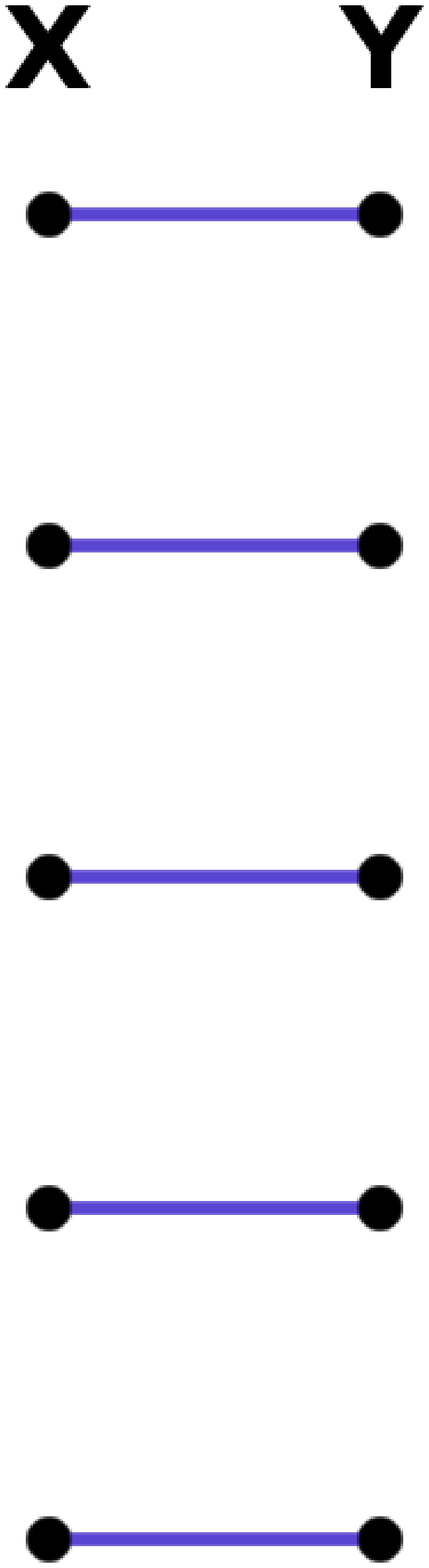}  
\includegraphics[scale=0.2]{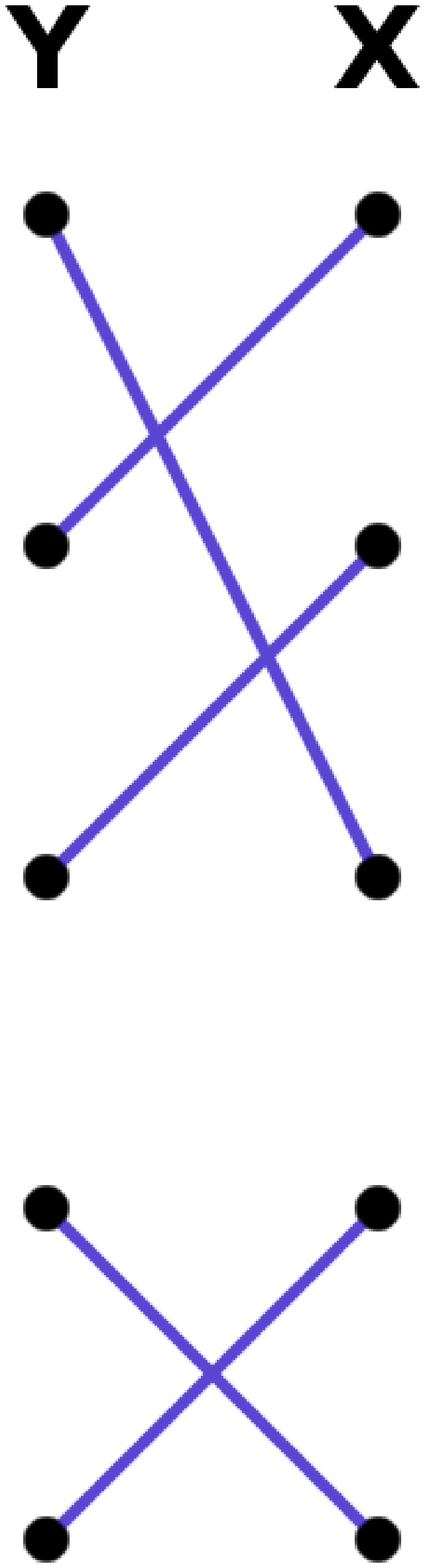}    
\includegraphics[scale=0.2]{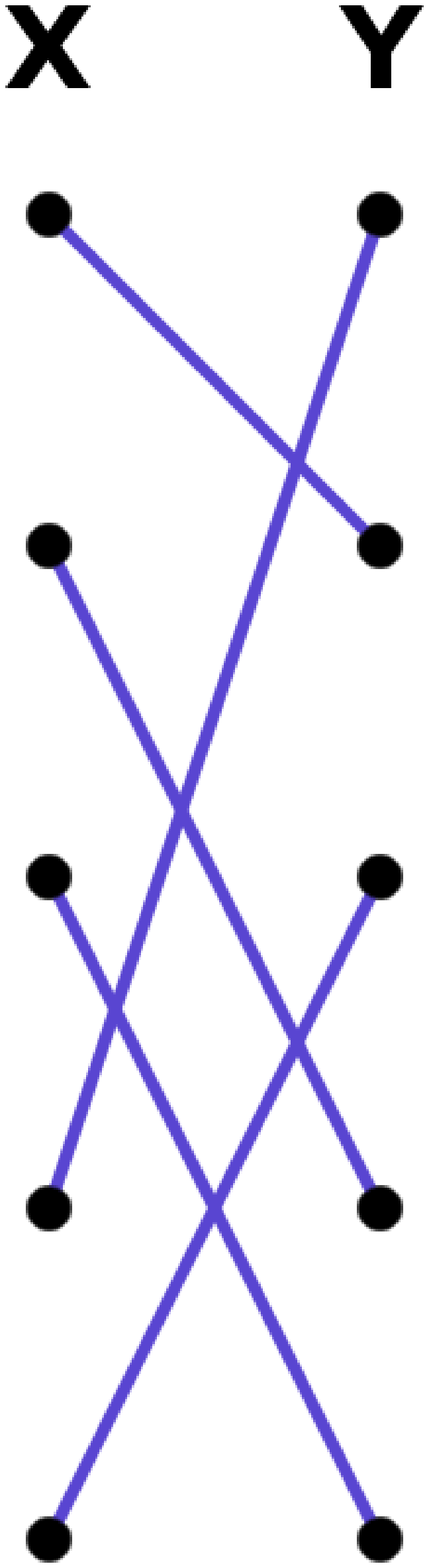} 
\includegraphics[scale=0.2]{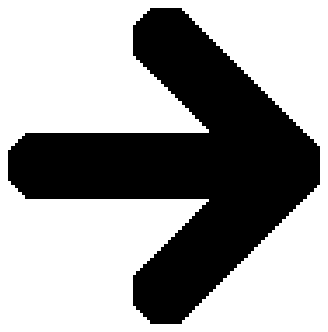} 
\includegraphics[scale=0.2]{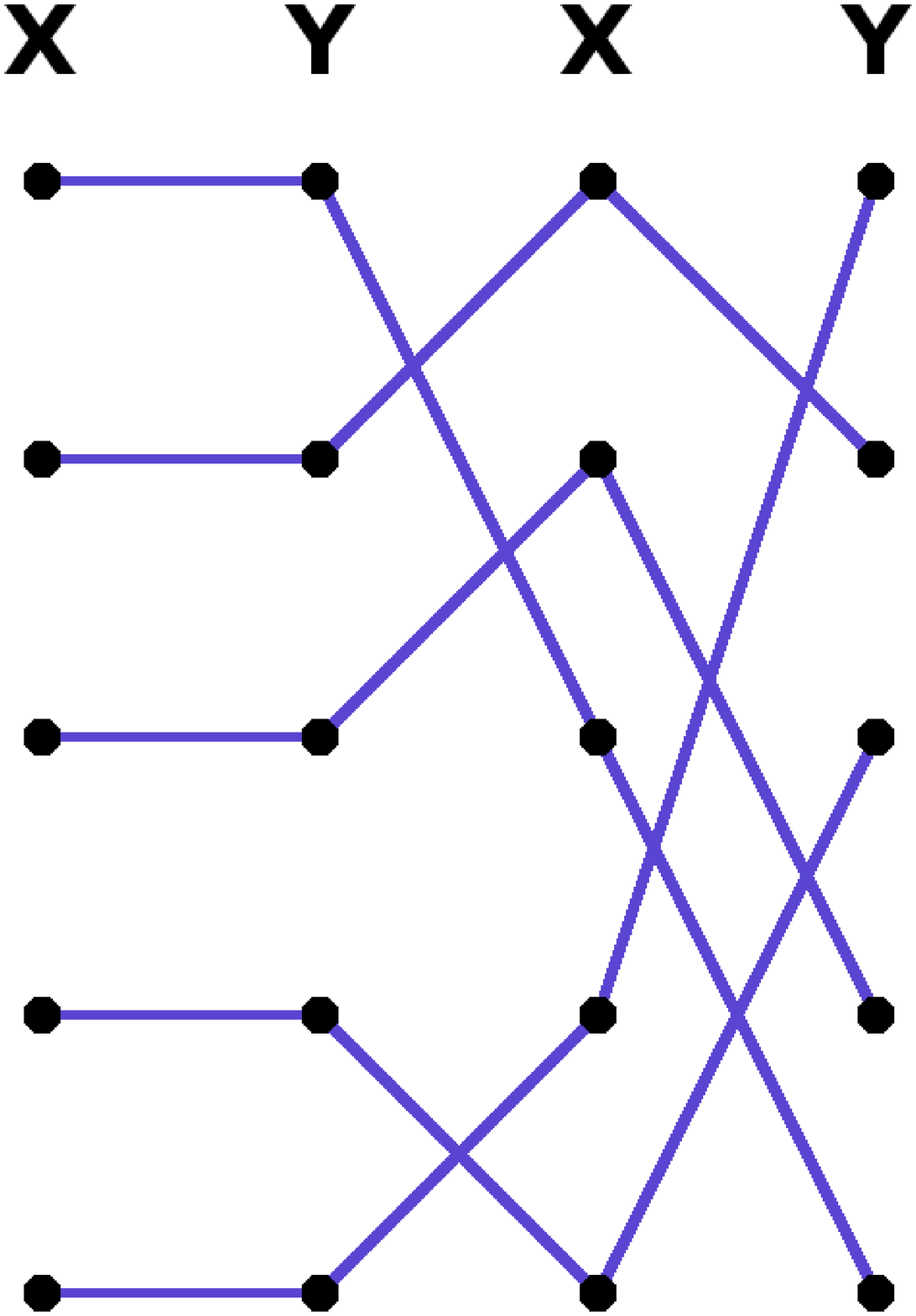}    
\caption{Concatenation of $3$ matchings, starting at $X$.}\label{f:Matching}
\end{figure}

Another way to decompose $\cc{Q}_n$ into paths, which will inform our method, is as follows. Since $\cc{Q}_n$ is $n$-regular and bipartite, it is a simple application of Hall's theorem that we can decompose the edge set into $n$ perfect matchings. Let $X$ be the set of even vertices in $\cc{Q}_n$ and $Y$ be the odd. If we take some perfect matchings $\cc{M}_1$, $ \cc{M}_2, \ldots,  \cc{M}_k$ then we can cover the edges in these matchings by $|X|$ walks of length $k$, one starting at each vertex in $X$. For example if the edge $(x_1,y_{i_1})$ is in $\cc{M}_1$ and the edge $(y_{i_1},x_{i_2})$ is in $\cc{M}_2$ and so on the we have that the walk starting at $x_1$ is $\{(x_1, y_{i_1}) ,(y_{i_1}, x_{i_2}) , \ldots (x_{i_{k-1}}, y_{i_k})\}$ if $k$ is odd, and $\{(x_1, y_{i_1}) ,(y_{i_1}, x_{i_2}) , \ldots (y_{i_{k-1}}, x_{i_k})\})$ if $k$ is even. We will use the notation $\cc{W}(\cc{M}_1, \cc{M}_2, \ldots , \cc{M}_k, X)$ to denote the set of walks formed by concatenating $\cc{M}_1$ to $\cc{M}_k$ in that order, starting at $X$, and similarly if we start at $Y$.  A pictorial representation of this process is presented in Figure \ref{f:Matching}.

Therefore, since as we noted before we can decompose $\cc{Q}_n$ into $n$ perfect matchings, we can use this method to decompose $\cc{Q}_n$ into walks of length $k$, for any $k\,|\,n$, by splitting the matchings into sets of size $k$ and concatenating them as above. If we are careful with the matchings we choose and the order we concatenate them in we can ensure that these walks are paths. For example if we take, for $1\leq i \leq n$, the matchings

\begin{equation}
\cc{M}_i = \{\left((q_1,q_2, \ldots, q_i , \ldots , q_n),(q_1,q_2, \ldots, q_i+1 , \ldots , q_n)\right) \,:\, (q_1,q_2, \ldots, q_i , \ldots , q_n) \in X\},
\end{equation}

where addition is performed modulo 2, we see that $\cc{W}(\cc{M}_1, \cc{M}_2, \ldots , \cc{M}_n, X)$ is exactly the antipodal paths of Lemma \ref{l:antipodal}.\\
\ \\
The main idea in the proof of Theorem \ref{t:main} is to first find a 'small' regular graph on $\cc{Q}_n$ which will interact nicely with paths we build up from matchings. It will be necessary to treat some small cases by hand and so we will take this opportunity to illustrate the ideas in the method with a small example. For example suppose we want to decompose $\cc{Q}_n$ into paths of length $2^2=4$.\\
\ \\
We first claim that, if we want to decompose $\cc{Q}_n$ into paths of length $4$, without loss of generality we can assume that $n \in [5,7]$. Indeed if $n \geq 9$ then $n-5 \geq 4$ and is even and so we have that $\cc{Q}_{n-5}$ can be decomposed into cycles of length $2^{n-5}$. Since $2^{n-5} > 4$ we can decompose each of these into paths of length $4$ and so $\cc{Q}_{n-5}$ can be decomposed into paths of length $4$. Therefore by Lemma \ref{l:split} it is sufficient to consider the cases where $n=5$ or $7$. We will just consider the case $n=5$ in this example. \\
\ \\
We view $\cc{Q}_5$ as $\cc{Q}_3 \times \cc{Q}_2$, that is for each $(p_1,p_2) \in \cc{Q}_2$ we look at the set of vertices $(q_1,q_2,q_3,q_4,q_5)$ such that $(q_4,q_5) = (p_1,p_2)$. The induced subgraph of $\cc{Q}_5$ on this set of vertices is $\cc{Q}_3$. We take a Hamiltonian cycle, $C$, on $\cc{Q}_3$ (it is a simple exercise to show that $\cc{Q}_n$ is Hamiltonian for all $n$) and take the union of these edges over all copies of $\cc{Q}_3$. That is, for each $(p_1,p_2) \in \cc{Q}_2$ we take the edge set of a copy of $C$ on the subgraph of $\cc{Q}_5$ restricted to the vertices $(q_1,q_2,q_3,q_4,q_5)$ such that $(q_4,q_5) = (p_1,p_2)$. We call the union of all these edges $G$, note that $G$ is a $2$-regular subgraph of $\cc{Q}_5$ which covers the vertices of $\cc{Q}_5$ with cycles of length $8$. Furthermore, since $\cc{Q}_3$ is $3$-regular, we have that $\cc{Q}_3 \setminus C$ is $1$-regular and bipartite, that is, it is a matching, $\cc{I}^*$. So the union over all copies of $\cc{Q}_3$ of $\cc{I}^*$, which we will denote by $\cc{I}$, is a matching on $\cc{Q}_5$. Since we have covered all the edges of each copy of $\cc{Q}_3$ with $\cc{I}$ and $G$, we have that the remaining edges of $\cc{Q}_5$ are just $\cc{M}_4$ and $\cc{M}_5$ from (2.1). We let $E(G) = E(G^0) \cup E(G^1)$, where $G^0$ is the restriction of $G$ to the vertices $(q_1,q_2,q_3,q_4,q_5)$ such that $q_4=0$ and similarly $G^1$ is the restriction to the vertices where $q_4=1$. This decomposition holds since all the edges of $G$ are contained within copies of $\cc{Q}_3$ inside $\cc{Q}_5$. We note that both $G^0$ and $G^1$ are $2$-regular, that is they cover their vertex set with cycles.\\
\ \\
We want to use $E(G^0)$ to extend $\cc{M}_4$ to paths of length $2$. We take each cycle in $G^0$ and arbitrarily give it an order by labelling the vertices. Given a cycle $C=\{(x_1,y_2),(y_2,x_3),(x_3,y_4), \ldots, (x_7,y_8),(y_8,x_1)\}$ in $G^0$ we look at the edges that are matched to $\{x_1,y_2,x_3,y_4,x_5,y_6,x_7,y_8\}$ in $\cc{M}_4$. Let us call them $y_1,x_2,y_3,x_4,y_5,x_6,y_7,x_8$ respectively, that is $(x_i,y_i) \in \cc{M}_4$ for $1\leq i \leq 8$. To each edge in the matching we adjoin the 'next' edge in the cycle, that is we form the set of paths

$$\big\{ \{(y_1,x_1), (x_1,y_2)\}, \{(x_2,y_2), (y_2,x_3)\}, \ldots, \{(x_8,y_8), (y_8,x_1)\} \big\}.$$

\begin{figure}[!ht]
\centering
\includegraphics[scale=0.2]{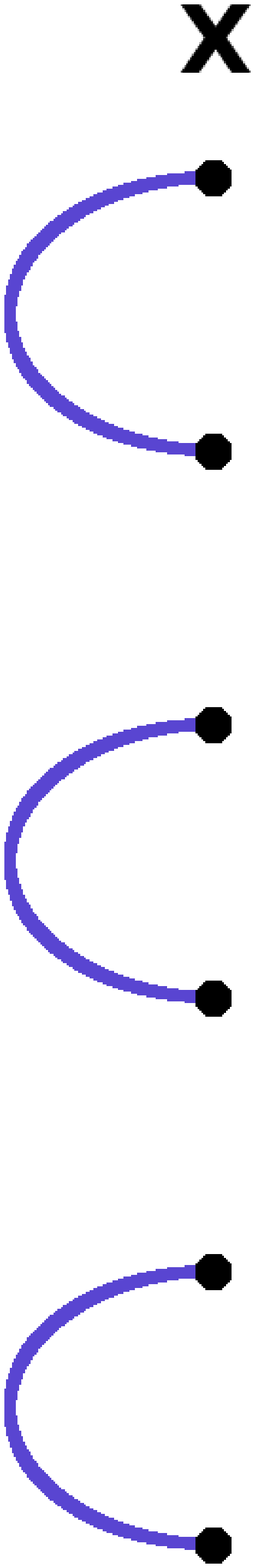}
\includegraphics[scale=0.2]{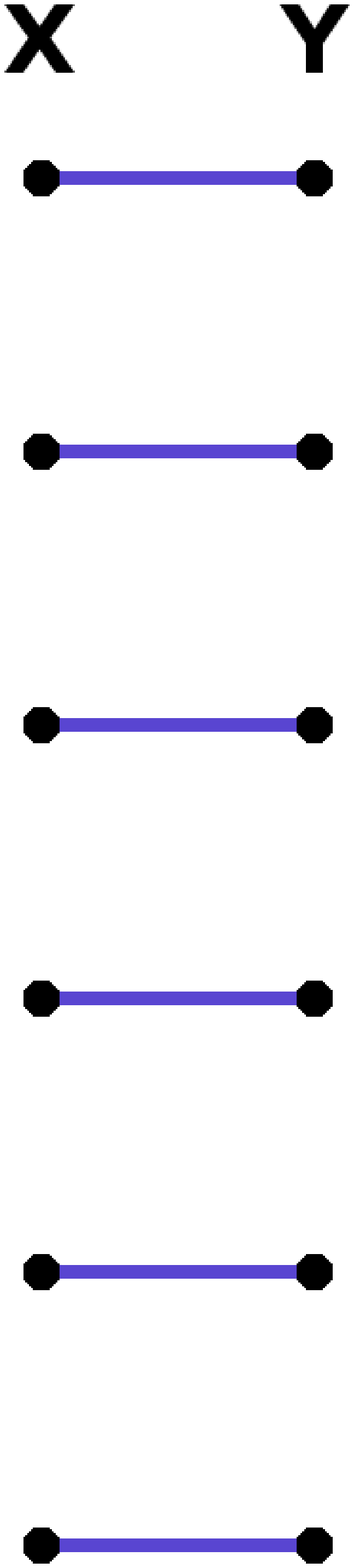}    
\caption{$\cc{P}_e$ and $\cc{I}$}\label{f:Join}
\end{figure}

We repeat this for every cycle in $G^0$, let us denote by $\cc{P}$ the union of these paths. Note that since $G^0$ is a graph on exactly half the vertices of $\cc{Q}_5$, and no edges of $\cc{M}_4$ are between vertices of $G^0$, we have that each edge of $\cc{M}_4$ is used in one of these paths, and also each edge of $G^0$. Since $\cc{Q}_5$ is bipartite we have that each of the paths in $\cc{P}$ is between two vertices from $X$, or two vertices from $Y$. In fact, moreover, each vertex in $X$ and $Y$ is an endpoint of exactly $1$ path, since $\cc{M}_4$ covered all the vertices of $\cc{Q}_5$ and $G^0$ was a union of cycles. Let us call the paths between even vertices \emph{even paths}, which we will denote by $\cc{P}_e$, and similarly the paths between odd vertices \emph{odd paths}, as $\cc{P}_o$. We will use these paths of length $2$ to join the edges in a matching into paths of length $4$. \\
\ \\
So we have that $E(\cc{Q}_5) = E(G) \cup E(\cc{I}) \cup E(\cc{M}_4) \cup E(\cc{M}_5) = E(G^1) \cup E(\cc{I}) \cup E(\cc{P}) \cup E(\cc{M}_5)$. Since $G^1$ consists of cycles of length $8$, it is simple to decompose it into paths of length $4$. We want to use $\cc{P}_e$ to connect the edges of $\cc{I}$ into paths of length $4$, and similarly $\cc{P}_o$ for $\cc{M}_5$. Each of the even paths is between two vertices in $X$, and each vertex in $X$ is used exactly once as an endpoint. Therefore we can form walks of length $4$ by adding to each path in $\cc{P}_e$ the two edges in $\cc{I}$ that are connected to it's endpoints, see Figure \ref{f:Join}. This will use each edge in $\cc{P}_e$ and $\cc{I}$. Similarly we use $\cc{P}_o$ to join the edges of $\cc{M}_5$ into walks of length $4$.\\
\ \\
We need to check that the walks we produce in this manner do not repeat any vertices. Since each path in $\cc{P}$ uses one edge from $\cc{M}_4$, given a walk $P = \{(x_1,x_2),(x_2,x_3),(x_3,x_4),(x_4,x_5)\}$ formed in this way, we have that without loss of generality either $x_1$ and $x_2$ are contained in the set of vertices with $q_4=0$ and $x_3,x_4$ and $x_5$ in the set of vertices where $q_4=1$, or $x_1,x_2$ and $x_3$ are contained in the set of vertices with $q_4=0$ and $x_4$ and $x_5$ in the set of vertices where $q_4=1$. Either way we have that $P$ consists of a walk of length 2 in one subcube and an edge in a disjoint subcube, joined together by an edge between these subcubes. Since the subcubes are disjoint, and each subcube is bipartite, no vertices are repeated and these walks are actually paths. Since we used all the edges of $\cc{P}$, $\cc{I}$ and $\cc{M}_5$ in this process, we have decomposed $\cc{Q}_5$ into paths of length $4$.\\
\ \\
For the general case our idea is similar, we will cover the vertices of a small subcube of $\cc{Q}_n$ with some cycles and then decompose the rest of the edges into two sorts of matchings, those contained inside copies of this subcube, like $\cc{I}$, and the rest of the form $\cc{M}_i$. We combine one of the $\cc{M}_i$ with some of the cycles from the subcube to form paths of length $2$, and join the rest of matchings into two sets of paths, one starting on $X$ and one starting on $Y$. We then use the paths of length $2$ as before to join the paths starting on $X$ pairwise, and similarly the paths starting on $Y$. If we have enough $\cc{M}_i$ compared to matchings from inside the subcube we can ensure that the walks we produce are actually paths, by making sure that in each walk we never use too many edges from inside the same copy of the small subcube.\\
\ \\

\begin{theorem}\label{t:pow}
Let $n$ be odd and $r$ such that $2^r < n$. Then $\cc{Q}_n$ can be decomposed into paths of length $2^r$.
\end{theorem}
\proof
Let us first suppose that $r$ is odd. By Lemma \ref{l:Ham} it is possible to decompose $\cc{Q}_{r+1}$ into Hamiltonian cycles, and so it is possible to decompose it into paths of length $2^r$, by splitting each cycle in half. Therefore by Lemma \ref{l:split}, it is sufficient to consider the case where $n= 2^r + l$ for some odd $1 \leq l \leq r$, note that $n=2^r + l \geq r+2$. We first build a subgraph on $\cc{Q}_n$ that is $l+1$-regular. For each $(p_1,p_2, \ldots,  p_{n-(r+1)}) \in \cc{Q}_{n-(r+1)}$ we consider the restriction of $\cc{Q}_n$ onto the set of vertices $(q_1,q_2, \ldots, q_n)$ such that $(q_{r+2}, q_{r+3}, \ldots,q_{n}) = (p_1,p_2, \ldots,  p_{n-(r+1)})$, this is isomorphic to $\cc{Q}_{r+1}$. By Lemma \ref{l:Ham} we can decompose $\cc{Q}_{r+1}$ into $\frac{r+1}{2}$ Hamiltonian cycles $C_1, \ldots, C_{\frac{r+1}{2}}$. We split each of the cycles $C_{\frac{l+3}{2}}, \ldots C_{\frac{r+1}{2}}$ into two matchings, so that we have decomposed the edge set of $\cc{Q}_{r+1}$ into $\frac{l+1}{2}$ cycles of length $2^{r+1}$, $C_1, \ldots, C_{\frac{l+1}{2}}$, and $r-l$ matchings, $\cc{I}^{*}_1, \ldots, \cc{I}^*_{r-l}$. For $1 \leq i \leq \frac{l+1}{2}$ we let $G_i$ be the graph formed by taking the union of the edge sets of a copy of $C_i$ on each copy of $\cc{Q}_{r+1}$. Similarly for all $1 \leq j \leq r-l$ we let $\cc{I}_j$ be the matching formed by taking a copy of $\cc{I}^*_j$ on each copy of $\cc{Q}_{r+1}$. We now have that

$$E(\cc{Q}_n) = \bigcup_{i=1}^{\frac{l+1}{2}} E(G_i) \cup \bigcup_{j=1}^{r-l} E(\cc{I}_j) \cup \bigcup_{t=r+2}^{n} E(\cc{M}_t).$$

As before we split $E(G_1)$ into $E(G_1^0) \cup E(G_1^1)$ where $G_1^0$ is the restriction of $G_1$ to the set of points $(q_1,q_2, \ldots ,q_n)$ such that $q_{r+2}=0$ and $G_1^1$ is the restriction of $G_1$ to the set of points $(q_1,q_2, \ldots ,q_n)$ such that $q_{r+2}=1$. Note that both $G_1^0$ and $G_1^1$ are $2$-regular graphs composed of a disjoint union of cycles of length $2^{r+1}$. We combine $G_1^0$ with $\cc{M}_{r+2}$ to form paths of length two as in the previous example. So, for every cycle $C=\{(x_1,y_2),(y_2,x_3),(x_3,y_4), \ldots, (x_{2^{r+1}-1},y_{2^{r+1}}),(y_{2^{r+1}},x_1)\}$ in $G_1^0$ we look at the edges matched to $\{x_1,y_2, \ldots , x_{2^{r+1}-1},y_{2^{r+1}}\}$ in $\cc{M}_{r+2}$. Let us call them $y_1,x_2, \ldots , y_{2^{r+1}-1},x_{2^{r+1}} $ respectively, that is $(x_i,y_i) \in \cc{M}_{r+2}$ for $1\leq i \leq 2^{r+1}$. To each edge in the matching we adjoin the 'next' edge in the cycle, that is we form the set of paths

$$\big\{ \{(y_1,x_1), (x_1,y_2)\}, \{(x_2,y_2), (y_2,x_3)\}, \ldots, \{(x_{2^{r+1}},y_{2^{r+1}}), (y_{2^{r+1}},x_1)\} \big\}.$$

We repeat this for every cycle in $G_1^0$, let us denote by $\cc{P}$ the union of these paths. As before we split $\cc{P}$ into a set $\cc{P}_e$ of even paths and a set $\cc{P}_o$ of odd sets, and note that every point of $X$ is an endpoint of exactly one path in $\cc{P}_e$, and similarly every point in $Y$ is an endpoint exactly one path in $\cc{P}_o$.\\
\ \\
We use the remaining matchings, $\cc{I}_1, \ldots \cc{I}_{r-l}$ and $\cc{M}_{r+3}, \ldots, \cc{M}_n$, to form two sets of walks, one starting at $X$ and one starting at $Y$, both of length $\frac{n-(l+2)}{2} = \frac{2^r-2}{2} = 2^{r-1}-1$. We want to order the matchings in such a way that these walks will be paths. For example if we took the set of walks $\cc{W}(\cc{M}_{r+3}, \cc{I}_1, \cc{M}_{r+4}, \cc{I}_2, \cc{M}_{r+5},\cc{I}_6, \ldots, X )$ alternating between using the $\cc{M}_i$ and the $\cc{I}_j$, at least until we run out of $\cc{I}_j$s, then the walks we form will actually be paths. Indeed, if we pick two vertices in the walk $x$ and $y$ which have an edge from $\cc{M}_i$, for some $i$, between them in the walk, then $x$ and $y$ do not agree in the $i$th co-ordinate. Therefore the only points that could be repeated in each walk are those joined by an edge in some $\cc{I}_j$, but clearly these are distinct, since $\cc{Q}_n$ has no loops. So we want to have at least as many $\cc{M}_i$s as we do $\cc{I}_j$s, that is we need that $n-(r+2) \geq r-l$. Since $n=2^r +l$ we need $2^r + 2l - 2 \geq 2r$ and since $l \geq 1$ it is sufficient that $2^r \geq 2r$, which holds for all odd $r$. So we form our two sets of paths in this way, one starting at $X$ and one starting at $Y$, and we use $\cc{P}_e$ to join the ones starting at $X$ and $\cc{P}_o$ to join the ones starting at $Y$ as before, into walks of length $2(2^{r-1}-1) +2 = 2^r$. Again it is a simple check that these walks are in fact paths. Let us consider one of the walks formed by combining $\cc{W}(\cc{M}_{r+3}, \cc{I}_1, \cc{M}_{r+4}, \cc{I}_2, \cc{M}_{r+5},\cc{I}_6, \ldots, X )$ and $\cc{P}_e$. It consists of two paths from $\cc{W}(\cc{M}_{r+3}, \cc{I}_1, \cc{M}_{r+4}, \cc{I}_2, \cc{M}_{r+5},\cc{I}_6, \ldots, X )$, joined together by a path of length two from $\cc{P}_e$. Since an edge of $\cc{M}_{r+2}$ was used in each path in $\cc{P}_e$ we have that the $2^{r-1}$ vertices in first path differ from the $2^{r-1}$ vertices in the second path in the $(r+2)$nd coordinate, and so they are all distinct. Finally the vertex in the middle of the path of length two differs from all of the vertices except it's immediate neighbours in the $(r+3)$rd coordinate, and since those three vertices were in a path in $\cc{P}_e$, it is distinct from those two as well. \\
\ \\
So, to conclude, we have decomposed $\cc{Q}_n$ into some graphs $G_1^1, G_2, G_3 \ldots G_{\frac{l+1}{2}}$ which are each a union of cycles of length $2^{r+1}$ and a collection of paths of length $2^r$, therefore $\cc{Q}_n$ can be decomposed into paths of length $2^r$.\\
\ \\
The case where $r$ is even is similar. Since we can decompose $\cc{Q}_{r+2}$ into Hamiltonian cycles it is sufficient to consider the case $n=2^r + l$ for some odd $1 \leq l \leq r+1$. We view $\cc{Q}_n$ as $\cc{Q}_{r+2} \times \cc{Q}_{n-(r+2)}$ and use the decomposition of $\cc{Q}_{r+2}$ into Hamiltonian cycles to split $\cc{Q}_n$ into $G_i$s, $\cc{I}_j$s and $\cc{M}_t$s as in the odd case.\\
\ \\
There are two small differences, firstly in order to make the paths of length $2$ we need that $\cc{M}_{r+3}$ exists. That is we need $n= 2^r + l \geq r+3$, but this holds for all even $r$, $1 \leq l \leq r+1$. The second difference comes when we want to check that we have at least as many $\cc{M}_i$s as $\cc{I}_j$s, since now we need that $n-(r+3) \geq r-l+1$, that is $2^r + 2l - 4 \geq 2r$. This holds for all $r \geq 4$, $1 \leq l \leq r+1$, and also for $r=2,l=3$. The only remaining case to check is therefore when $r=2$ and $l=1$, that is, we need to demonstrate a decomposition of $\cc{Q}_5$ into paths of length $4$, which we did in the preceding example.
\qed\\
\ \\

\emph{Proof of Theorem \ref{t:main}} 
Given $k\,|\,n2^{n-1}$ we have that $k=t2^r$ for some odd $t\,|\,n$. Since $k=t2^r \leq n$ we have that $2^r \leq \frac{n}{t}$, and so by Theorem \ref{t:pow} $\cc{Q}_{\frac{n}{t}}$ can be decomposed into paths of length $2^r$. Therefore by by Lemma \ref{l:decomp} $\cc{Q}_n$ can be decomposed into paths of length $k$. \qed
\ \\

\section{Even $n$}
The case where $n$ is even seems different. For example in the odd case the problem seems just as difficult if we ask for walks instead of paths. However for even $n$, since every vertex has even degree, $\cc{Q}_n$ has an Eulerian cycle, and so it is possible to decompose $\cc{Q}_n$ into walks of length $k$ for every $k\,|\,n2^{n-1}$. If we want to decompose $\cc{Q}_n$ into {\bf paths} of length $k$ we still need that $k\,|\,n2^{n-1}$, but we no longer require that $k \leq n$. For example by Lemma \ref{l:Ham} we can decompose $\cc{Q}_n$ into paths of length $2^{n-1}$, so a more natural condition would seem to be $k < 2^n$, since no path can be longer than $|\cc{Q}_n|$. The methods of Section 2 prove some results towards this, for example since Lemma \ref{l:antipodal} and Lemma \ref{l:decomp} hold for general $n$ we know that if $t\,|\,n$ with $t$ odd, then we can decompose $\cc{Q}_n$ into paths of length $t2^{\frac{n}{t}-1}$. However it may be possible to decompose it into paths of length $t2^{n-\lceil \log{t} \rceil}$. We conjecture:

\begin{conjecture}
Let $n$ be even and $k$ such that $k \,|\, n2^{n-1}$ and $k < 2^{n}$. Then $\cc{Q}_n$ can be decomposed into paths of length $k$.
\end{conjecture}

\bibliography{Cube}
\bibliographystyle{plain}
\end{document}